\documentclass{amsart}
\usepackage{latexsym,amsxtra,ifthen}
\usepackage{amsfonts}
\usepackage{verbatim}
\usepackage{amsmath}
\usepackage{amsthm}
\usepackage{amssymb}

\numberwithin{equation}{section}

\theoremstyle{plain}
\newtheorem{theorem}{Theorem}[section]

\newtheorem{prop}[theorem]{Proposition}
\newtheorem{lemma}[theorem]{Lemma}
\newtheorem{cor}[theorem]{Corollary}

\theoremstyle{definition}

\newtheorem{nota}[theorem]{Remark}

\newcommand{\N}{\ensuremath{\mathbb N}}

\newcommand{\Z}{\ensuremath{\mathbb Z}}

\def\t{\widetilde}

\def\deg{\mbox{deg\,}}
\def\dim{\mbox{dim\,}}

\def\ext{\mbox{Ext\,}}
\def\gd{\mbox{gl. dim\,}}

\def\gr{\mbox{\bf gr\,}}

\def\id{\mbox{id\,}}

\def\tor{\mbox{Tor\,}}

\begin{document}

\title{Koszul algebras associated to graphs%
}

\author{Dmitri Piontkovski}

      \address{Department of High Mathematics for Economics, 
Myasnitskaya str. 20, State University `Higher School of Economics', Moscow 101990, Russia 
}
\thanks{Partially
supported by the grant 02-01-00468 of the Russian Basic Research Foundation}

\email{piont@mccme.ru}



\date{\today}

\begin{abstract}
Quadratic algebras associated to graphs have been introduced by I. Gelfand, S. Gelfand, and Retakh
in connection with decompositions of noncommutative polynomials. 
Here we show that, for each graph with rare triangular subgraphs, 
the corresponding quadratic algebra is a Koszul domain with global dimension equal to 
the number of vertices of the graph. 
\end{abstract}

\maketitle

\section{Introduction}


The origin of the algebras considered here are the works of I. Gelfand, S. Gelfand, Retakh, Wilson, and others
on decompositions of noncommutative polynomials and noncommutative symmetric
functions, see~\cite{qdet} and references therein. Here we sketch some their definitions and results.

Let $C$ be a division algebra over a field $k$, and let $P(x) \in C[x]$
be a polynomial of degree $n$, where $x$ ia a central variable. 
An element $t \in C$ is called {\it pseudo-root}
of the polynomial $P(x)$ if there is a decomposition 
$$
P(x) = L(x) (x-t ) R(x),
$$
where $L(x), R(x) \in C[x]$ are polynomials. 
For a generic monic polynomial $P(x)$, there are at least $n2^{n-1}$ different pseudo-roots,
but these pseudo-roots are connected by some linear and quadratic relations. 
So, it is natural to consider an {\it universal quadratic algebra of pseudo-roots} 
$Q_n$, which is generated by these generic pseudo-roots 
satisfiing these relations. 

The algebra $Q_n $ has been studied in several papers~\cite{grw, 5, sw, pi1}. 
It is a homogeneous quadratic Koszul algebra with $2^n-1$ generators
and is very far from being commutative.
These generators $u_S$ may be naturally indexed by the nonempty subsets of the $n$-element set
$S \subset \{1, \dots , n \} $.  If we consider these subsets $S$ as faces of the $n-1$-dimensional 
simplex $\Delta$, 
it is natural consider the ideals of $I_\Gamma$ in $Q_n$ generated by all simplicial subcomplexes 
$\Gamma$ (that is, by the variables $u_S$ with $S\in \Gamma$) 
of the simplex $\Delta$. The quotient algebras $Q(\Gamma) = Q_n/I_\Gamma$ 
would lead to ``noncommutative combinatorial topology''~\cite{qdet}. 
We need to ``glue'' such algebras corresponding to different simplicial complexes 
in order to approximate the noncommutative structure of the algebra $Q_n$.  
To do this,  we need to know, at first, the Hilbert series of the algebras 
$Q(\Gamma)$. Next, the algebra $Q_n$ is Koszul and has global dimension $n$~\cite{sw, pi1}.
Therefore, to approximate it by some algebras $Q(\Gamma)$, it is important to answer 
a question stated
by Vladimir Retakh: For which simplicial complexes $\Gamma$ the algebra $Q(\Gamma)$ is Koszul?    
We hope that all such algebras are Koszul of global dimension $n$. 


The generators and relations of the algebras $Q(\Gamma)$ has been described in~\cite{CW}, where 
these algebras have been introduced. In the trivial case on 0-dimensional 
simplicial complex $\Gamma$ (a set of $n$ points), 
the algebra $Q(\Gamma)$ is an algebra of commutative polynomials in $n$
variables.  
Special attention was given to the first non-trivial case of 
one-dimensional simplicial complex $\Gamma$ (i.~e., the case of a graph): 
correspondent  algebras seems to be the "most commutative" among all algebras $Q(\Gamma)$.
In particular, in the case of one-dimensional $\Gamma$  a simple set of relations has been described in~\cite{CW}.

The second step has been made by Nacin~\cite{nacin1, nacin2, nacin3}.  
He proved that the algebra $Q(\Gamma)$ is Koszul 
for graphs of the following types: ``line'' (i.~e., a graph with $n$ vertices $1, \dots , n$ 
and the edges $(12),(23),\dots , (n-1,n)$), ``star'' (vertices are $(12),(13),\dots , (1n)$), ``triangle''
(a complete graph with vertices). It is shown that, in these 3 cases, the global dimension of the algebra 
$Q(\Gamma)$ is equal to the number of vertices, and Hilbert series of the algebra and its 
quadratic dual algebra are calculated. In particular, the following interesting fact 
is found in the  cases of lines and stars: If $Q(\Gamma)^!(z) = p_0 +p_1 z + \dots + p_nz^n$ 
is the Hilbert series of the quadratic dual algebra $Q(\Gamma)^!$ (i.~e., $p_i = \dim Q(\Gamma)^!_i$), 
then the ``paliandromic''
equalities $p_{i} = p_{n-i}$ hold for all $i \le n$. However, in the case of triangle 
this paliandromic property fails.

In this paper, we try to make a next step in this direction. 
We consider a more general class of graphs, that is, the class of graphs 
with no overlapping triangles, i.~e., our graphs do not contain  
subgraphs isomorphic to the ``simple butterfly'' $\{ (12),(13),(23),(14),(15),(45) \}$
and the ``square with a diagonal'' $\{ (12),(23),(34),(14), (13) \}$.
In particular, this class includes all trees and all cyclic graphs. 

Our main result is the following.

\begin{theorem}[Theorem~\ref{Koszulity}]
Suppose that a graph $\Gamma$ with $n$ vertices
 does not contain a pair of overlapping triangles. 
Then the algebra $Q(\Gamma)$ is a Koszul domain of global dimension $n$.
\end{theorem} 

Also, we give explicit formulae for Hilbert series of algebras $Q(\Gamma)$ 
and their dual algebras $Q(\Gamma)^!$, see Corollary~\ref{Hilb} and Proposition~\ref{basis}
below. In particular, we explain the phenomenon of ``paliandroms'' found in~\cite{nacin2, nacin3}
for some graphs.  

\begin{cor}[Corollary~\ref{dual_coef}]
Let $p_t = \dim Q(\Gamma)^!_t$ denotes the Hilbert function of the quadratic
algebra $Q(\Gamma)^!$  dual to $Q(\Gamma)$. If $\Gamma$ does not contain two overlapping triangles, 
then $p_t =0$ for $t>n$ and 
$p_t \ge p_{n-t}$ for all $t \le n/2$.  The equalities $p_t = p_{n-t}$  hold for all  $t \le n/2$
if and only if the graph $\Gamma$ is triangle free.
\end{cor}

Let us say a few words about other probable properties of algebras $Q(\Gamma)$.
If the graph $\Gamma$ has at least one edge,  then the algebra $Q(\Gamma)$
has exponential growth (since its quotient algebra $Q_2$ has exponential growth), 
hence  $Q(\Gamma)$ is not Noetherian~\cite{SZ}. However, we conjecture that the algebras $Q_n$ and $Q(\Gamma)$
are (graded) coherent, that is, 
the kernel of every (homogeneous) map $M\to N$ of free
finitely generated $Q(\Gamma)$--modules is finitely generated. At least, the algebra $Q_2$
is coherent by~\cite{pi2}. If it is the case, the Koszul duality gives equivalences of 
several categories of finitely presented $Q(\Gamma)$--modules 
and their derived categories with the respective categories of modules over the finitely-dimensional 
algebra 
$Q(\Gamma)^!$~\cite{bgs, mvs}.

 \subsection*{Acknowledgement}
 
 I am grateful to Vladimir Retakh for helpful discussions. 
  I am also grateful to MPIM Bonn
 for their hospitality during preparation of this paper.

\section{Background}

\subsection{Quadratic algebras corresponding to graphs}

\label{def_QG}

Let $\Gamma = (W,E)$ be a graph with $n$ vertices $1,
\dots, n$, that is, 1--dimensional simplicial complex. A
construction in~\cite{CW} assigns to $\Gamma$ a quadratic
associative algebra $Q(\Gamma)$ with generators $\{ u_i |
i= 1 \dots n \} \cup \{ u_{ij} | 1 \le j i \le n, (ij) \in
E \}$ (where a sub-index ``$ij$'' denotes an unordered
pair) and relations for all pairwise different indexes
$i,j,k,l$
$$
    R(ij) := [u_i, u_j] - u_{ij} (u_i - u_j), \mbox{ where } i \ne j,
$$
$$
    R(ijk) := [u_i, u_{jk}] +[u_{ik},u_j]+ [u_{ik},u_{jk}]
      - u_{ij} (u_{ik} - u_{jk}),
$$
$$
    R(ijkl) : = [u_{ij}, u_{kl}]
$$
 (where $u_{pq} = 0$ for $(pq) \notin E$). 

\subsection{Koszulity}

Recall that a graded connected degree-one generated algebra $A$
over a field $k$ (i.~e. $A = k \oplus A_1 \oplus A_2 \oplus \dots$ is generated by the 
graded component $A_1$) is called {\it Koszul}
if the trivial right $A$-module $k_A$ has a linear free resolution, that is, 
$\tor^A_i(k,k)_j = 0$ for all $i,j$ but $i=j$ (graded the homologies of $A$ are 
concentrated in the diagonal).
We will use the following criterion for Koszul algebras~\cite{pp}.

Let $A$  be a quadratic algebra.
Assume that $A$ is filtered by a ${\N}$--graded ordered semigroup $S$
such that the filtration on the  graded component  $A_n$ is induced by the filtration on
the space of generators $V = A_1$.
Then $A$ has also an ${\N}$--filtration induced by the grading of  $S$.
The associated graded algebra
$\gr A$
may in general be non-quadratic;
let $R \subset A_1 \otimes A_1$ be a set of quadratic relations of
the algebra $\gr A$ (that is, the degree two graded component of the ideal of relations) 
and let $B = T(A_1)/\id (R)$ be its quadratic part.

\begin{theorem}[\cite{pp}, Theorem~7.1 in Ch.~4]
\label{pp-th}
 Assume that, in the notation above,

(i) the algebra $\gr A$ has no nontrivial relations in degree 3

and

(ii) the algebra $B$ is Koszul.

Then $\gr A = B$ and $A$ is Koszul.
\end{theorem}

The condition $(i)$ here is equivalent to the equality $\dim A_3 = \dim B_3$, because
$\dim A_3 = \dim \gr A_3$.

A classical particular case of this theorem is that
every algebra with quadratic Groebner basis is Koszul~\cite{pri}
 (since every monomial quadratic algebra $B$ is Koszul).
Another standard consequence
is that every quotient of  a polynomial ring
or an  external algebra by an ideal with quadratic Groebner basis
is Koszul (because the Koszulness of a quotient of an external algebra by
a set of quadratic monomials follows from the results of~\cite{fro}).
We will use this corollary in the case of external algebras.

For the background on Groebner basis in the free associative
algebra, we refer the reader to~\cite{ufn}.

\section{Associated quadratic algebra to $Q(\Gamma)$}

Let us introduce a filtration on the algebra $A = Q(\Gamma)$
by setting new degrees of generators as $|u_i| = 1$ and $|u_{ij}| = 0$.
Then the associated quadratic algebra $B = B(\Gamma)$
has the same generators as $A$ and the relations
$$
    S(ij) := [u_i, u_j],
$$
$$
    S(ijk) := [u_i, u_{jk}] +[u_{ik},u_j],
$$
$$
    S(ijkl) := R(ijkl) = [u_{ij}, u_{kl}]
$$ (where again $u_{pq} = 0$ for $(pq) \notin E$, and the indexes
$i,j,k,l$ are pairwise different). Because these relations are
linear combinations of commutators of generators, the algebra $B$
is a universal enveloping algebra of the quadratic Lie algebra $g
= g(\Gamma)$ given presented by the same generators and relations.

To use Theorem~\ref{pp-th}, we need

\begin{prop}
\label{dim3}
 For every graph $\Gamma$, we have $\dim Q(\Gamma)_3 = \dim
B(\Gamma)_3$.
\end{prop}

\begin{proof}
Let $u$ be the minimal set of degree-one homogeneous generators of
$A$ (and $B$) listed above. Let $G = \{ 0,1 \}$ be the
two--elements multiplicative semigroup. 
Let us introduce  a
${G}^n$-grading on the free algebra $T(u)$ by assigning $\deg' u_i
= {\bf 1}_i := (0, \dots, 1, \dots, 0)$ (the unit in the $i$-th
place) and $\deg' u_{ij} : = {\bf 1}_i + {\bf 1}_j$.  
Then the relations
of both algebras $A$ and $B$ becomes homogeneous w.~r.~t. this new
grading, so, the new ${G}^n$ -rading is correct, and the algebras 
become ${\Z}\times {G}^n$-graded. 

 In order to check the equality 
$\dim Q(\Gamma)_3 = \dim
B(\Gamma)_3$, we have to check the equalities 
for dimensions of graded components $\dim Q(\Gamma)_{3,g} = \dim
B(\Gamma)_{3,g}$ for every $g \in {G}^n$.
For $g = (g_1, \dots, g_n)\in {G}^n$, let $|g|$ denotes the 
number of units among $g_i$. Because for every generator $u_*$
we have $|\deg' u_*| \le 2$, hence $Q(\Gamma)_{1,g} = 
B(\Gamma)_{1,g} = 0$ if $|g| > 2$, therefore, 
$ Q(\Gamma)_{3,g} = 
B(\Gamma)_{3,g} = 0$ for $|g| > 6$. 
So, it remains to check the equalities $\dim Q(\Gamma)_{3,g} = \dim
B(\Gamma)_{3,g}$ for all $g$ with $|g|\le 6$.
This means that there are some subindexes  $i_1, \dots, i_s$, where $s 
\le 6$, such that all elements in $Q(\Gamma)_{3,g} $ and
$ B(\Gamma)_{3,g}$  lies in the span of monomials on some generators $u_*$
depending on these subindexes only. Hence we may
ignore all the relations of algebras $A$ and $B$ which depend on
other subindexes, that is, we work in the subgraph of $\Gamma$ with
vertices $i_1, \dots, i_6$. This means that {\it it is sufficient
to prove proposition~\ref{dim3} for $n \le 6$}.

Moreover, if  $|g| = 6$, then all the monomials in the decomposition of 
every element $f\in Q(\Gamma)_{3,g}$ (or $f\in B(\Gamma)_{3,g}$)
 must
have the form $\mu = u_{ij} u_{kl} u_{mn}$, where all indexes $i, \dots , n$
are pairwise different.
Let $\Xi$ 
be the span of monomials of  that type in the free algebra $T(A_1)$ generated by $u$.
The intersections of $\Xi$  with the ideal of relations 
of the algebras $Q(\Gamma)$ 
and $B(\Gamma)$ are the same as its intersection
with the ideal  $I_\Xi$ generated by the relations of the type $R(ijkl)$, 
because $R(ijkl) \in \Xi$, while for every other relation $R_*$ or $S_*$ listed above,
  any monomial in the decomposition of this noncommutative  polynomial 
  is not  a submonomial of a monomial of  the form $\mu$.
  Therefore, 
  the vector spaces $Q(\Gamma)_{3,g}$ are $B(\Gamma)_{3,g}$ are isomorphic  to the 
  correspondent  graded component of the algebra $T(A_1)/I_\Xi$, hence they are isomorphic 
  to each other.

Now, it remains to consider the case $|g| \le 5$.
By the same arguments as above, it sufficient 
to show the following

\begin{lemma}
\label{dim3n5} If $n \le 5$, then the statement of Proposition~\ref{dim3} is true.
\end{lemma}

\begin{proof}[Proof of Lemma~\ref{dim3n5}]
Being a statement about dimensions of some subspaces in
finite--dimensional vector spaces, this Lemma is proved by a
direct computer calculation. 
There are 33 pairwise non-isomorphical graphs with $\le 5$ vertices.
For all these cases, suitable dimensions of the third components of algebras 
$Q(\Gamma)$ 
and $B(\Gamma)$ have been calculated via computations of Hilbert series up to the 
3-rd power. For this purpose, I used a software GRAAL 
(developed by Alexei Kondratiev).
In each of these 33 cases, the calculation shows that 
the Hilbert series of the algebras $Q(\Gamma)$ 
and $B(\Gamma)$  are equal to each other at least up to $o(z^4)$;
in particular,   $ \dim Q(\Gamma)_3 = \dim
B(\Gamma)_3$.
\end{proof} 

So, Proposition~\ref{dim3} is proved completely.
\end{proof} 

In the view of Theorem~\ref{pp-th}, we have

\begin{cor}
\label{BA}
 If, for some $\Gamma$, the algebra $B(\Gamma)$ is Koszul,
then the algebra $Q(\Gamma)$ is Koszul too, and both algebras are domains.
\end{cor}

\begin{proof}
The Koszulity is a part of  Theorem~\ref{pp-th}. 
It remains to show that the algebras $Q = Q(\Gamma)$ and $B =B(\Gamma)$ have no zerodivisors.
Indeed, $B$ is an universal enaveloping algebra of some graded Lie superalgebra
(because the relations of the algebra $B$ are linear combinations of commutators of some variables), 
hence it is a domain. By Theorem~\ref{pp-th}, $B = \gr A$, so,  $A$ is a domain too.
\end{proof}

\section{Graphs without overlapping triangles}



\begin{theorem}
\label{Koszulity}
Suppose that there do not exist two triangles with common vertex
in $\Gamma$.
 Then both algebras $A = Q(\Gamma)$ and $B = B(\Gamma)$ are Koszul domains 
 and have global dimension $n$. 
\end{theorem}

\begin{proof}
The dual quadratic algebra $B^!$ is a quotient
of an external algebra $\Lambda V$
(where $V$ is a span of indeterminates
$\{ e_i | i= 1 \dots n \} \cup
\{ e_{ij} | 1 \le j  i \le n, (ij) \in E \}$)
by an ideal generated by the following relations:
\begin{equation}
\label{RelB!}
  \begin{array}{c}
    T(ijk) := e_i e_{jk} + e_j e_{ki} + e_k e_{ij},
    \mbox{ where } (ij), (jk), (ij) \in E, \\
    U(ik) := e_i e_{ik}, \\
    W(ijk) := e_{ik} e_{jk}.
   \end{array}
\end{equation}

 We claim that these relations (with $i < j$) form a
Groebner basis of the ideal $I$ in $\Lambda V$  generated
by them w.~r.~t. any degree--lexicographical order with
$e_n >\dots e_1 > e_{ij}$.

We use standard Buchberger criterion. 

A nontrivial s--polynomials may occur only involving some
non--monomial relation $T(ijk)$, where $(jk), (ik), (ij) \in E$
with $i>j, i>k$. Its leading term is $\t {T(ijk)} = e_i
e_{jk}$. The following
 s--polynomials  are obviously reduced  to zero
by
$U(ik)$ and $W(ijk)$ (in these presentations, we omit the monomials ovbiously
reduced to zero):
$e_i T(ijk) = e_ie_j e_{ik} + e_ie_k e_{ij}$,
 $T(ijk) e_{jk} = e_j e_{ik}e_{jk} + e_k e_{ij}e_{jk}$,
 $T(ijk) e_{im}
+U(im)e_{jk} = e_j e_{ik}e_{jm} +e_k e_{ij}e_{jm}$,
 $T(ijk) e_{k} + e_iU(kj) = e_{ik} e_j e_{k} $,
 and $T(ijk)e_{sk} +e_i W(sjk) = e_j e_{ik}e_{sk} +e_k e_{ij} e_{sk}$.

It remains to consider the intersections between leading monomials
of $T(ijk)$ and some other $T(pqr)$. If $T(pqr)$ exists, then the
triangle $(pqr) = \{ (pq),(qr),(pr) \}$ exists in $E$.
If there is an overlap between the leading terms of $T(ijk)$ and
$T(pqr)$, the triangles $(ijk)$ and $(pqr)$ must have at least one
common vertex. But it is impossible, so that there is no other
s-polynomials.

So, the algebra $B^!$ has quadratic Groebner basis as a
quotient of $\Lambda V$. Hence it is Koszul, so that $B$ is
Koszul too. By Corollary~\ref{BA}, both algebras $A$ and $B$ are Koszul domains.
Because $\gd A = \max \{ d | A^!_d \ne 0 \}$ and $A^!(z) =
A(-z)^{-1} = B(-z)^{-1} =B^!(z)$, we have $\gd A = \gd B =
\max \{ d | B^!_d \ne 0 \}$. To calculate this global
dimension, we give a description of linear basis of the
algebra $B^!$.

By the definition of the Groebner basis, there is a  linear basis
of $B^!$ consisting of all monomials in $\Lambda V$ which are not
divisible by the leading terms of the Groebner basis. These
leading terms are the following:
\begin{equation}
\label{GroB!}
  \begin{array}{c}
       e_j e_{jk}; \\
       e_i e_{jk}, \mbox{ where $i>j>k$  and }
              (ij), (ik), (ij) \in E; \\
        e_{ij} e_{jk}, \mbox{ where } i >k.
   \end{array}
\end{equation}
 In particular, the only element of degree $n$ in the
basis of $B^!$ is $e_1 \dots e_n$, and there are no
elements of higher degree.
\end{proof}

Notice that we have obtain a linear basis of the algebra $B(\Gamma)^!$. Since the algebra $Q(\Gamma)^!$
is an associated graded algebra to a filtration on it~\cite[Ch.~4, Corollary~7.3]{pp}. 
So, we get

\begin{cor}
\label{basis-QG}
The monomials of the form 
$$ e_{i_1 i_2}, \dots ,e_{i_{2p-1} i_{2p}}  e_{j_1} \dots e_{j_{q}},
$$ where 
where $i_1 > i_3> \dots > i_{2p-1}$, ${j_1} > \dots > {j_{q}}$,  
$i_1 > i_2, \dots, {i_{2p-1} > i_{2p}}$, all indexes $i_1, \dots, i_{2p}, j_1, \dots , j_q$
are pairwise different, and for every $s \le p, t \le q$ either 
$j_t < i_{2s-1} $ or $(j_t, i_{2s-1}) \notin E$ 
or  $(j_t, i_{2s}) \notin E$,
form a linear basis of the algebra $Q(\Gamma)^!$. 
\end{cor}

\begin{cor}
\label{Hilb}
Suppose that there do not exist two triangles with common vertex
in $\Gamma$.
Then the Hilbert series of the algebra $Q(\Gamma)$ is 
$$Q(\Gamma)(z) = p_{\Gamma}(-z)^{-1}, $$
where $p_{\Gamma}(z) = Q(\Gamma)^!(z) = B(\Gamma)^!(z)$
is a polynomial of degree $n$ with positive integer coefficients.
\end{cor}

To give an implicit  formula for this polynomial $p_\Gamma(z) = p_0 + p_1 z + \dots +p_n z^n$,
 let us introduce some notations. 
Let $L_p$ be the set of all $p$-elements subsets 
$W=\{ e_{i_1 i_2}, \dots ,e_{i_{2p-1} i_{2p}}  \}$ (where $i_{2j-1} > i_{2j}$ for all $j$) 
such that 
$i_s \ne i_t$ for all $s,t= 1 \dots 2p$.
Given such a subset $W$, let $R_W$ denote the set of indexes $m \in [1..n]$ such that 
$m \ne i_s$ 
 and for every $j = 1 \dots p$ either $m < i_{2j-1} $ or $(m, i_{2j-1}) \notin E$ 
or  $(m, i_{2j}) \notin E$ 
(by the other words, we avoid all triangles $\{ m, i_{2j-1}, i_{2j} \}$
with $m> i_{2j-1}> i_{2j}$). 
Let
 $l_p$ and $r_W$ denote the cardinalities of the sets $L_p$ and $R_W$.


\begin{prop}
\label{basis}
Suppose that there do not exist two triangles with common vertex
in $\Gamma$. Then the Hilbert series $p(z) = Q(\Gamma)^!(z) = B(\Gamma)^!(z)$
is calculated by the formula 
$$
         p(z) = \sum_{p=0}^{[n/2]} z^{p}  \sum_{W\in L_p} (1+z)^{r_W}.
$$
In particular, if there is no triangle in  $\Gamma$, then 
$$     
p(z) = \sum_{p=0}^{[n/2]} l_p  z^p (1+z)^{n-2p}.
$$
\end{prop}

\begin{proof}
The coefficient $p_t$ is equal to the number of elements in the $t$-th
graded component of the linear basis of the algebra  $B^!$ described above. 
It follows from the description~(\ref{GroB!}) that this graded component 
consists of the square free commutative monomials of the form 
$  e_{i_1 i_2}, \dots ,e_{i_{2p-1} i_{2p}}  e_{j_1} \dots e_{j_{t-p}}$, \
where  $W:=\{ e_{i_1 i_2}, \dots ,e_{i_{2p-1} i_{2p}}  \} \in L_p$ 
and $j_s \notin R_W$
for all $s = 1 \dots t-p$.
Given such $W$ and $t>p$, we obtain exactly $\binom{r_W}{t-p}$ monomials of this type.
Now, the desired equalities immediately follow from the binomial decompositions.  
\end{proof}

In particular, we have 

\begin{cor}
\label{dual_coef}
If $\Gamma$ does not contain two overlapping triangles, we have
$p_t \ge p_{n-t}$ for all $t < n/2$.  The equalities $p_t = p_{n-t}$  hold for all  $t < n/2$
if and only if there is no triangle in $\Gamma$.
\end{cor}

\begin{proof}
The first inequality holds (for every $t$) for the coefficients of $z^t$ and $z^{n-t}$ 
in every summand of the first 
sum in Proposition~\ref{basis} (because $r_W \le n-2p$), 
hence it holds also for the whole sum. The equality means that 
$r_W = n-2p$ for every $W\in L_p$, that is, there is no triangle in $\Gamma$.
\end{proof}

\begin{nota}
In the view of `paliandromic' equalities $p_t = p_{n-t}$ 
for the dimensions of graded components of algebras $Q(\Gamma)^!$ and $B^!$ 
for triangle free ghaphs, one could conjecture that these algebras are Frobenius;
in this case, the algebras  $Q(\Gamma)$ and  $B(\Gamma)$  
are Gorenstein of finite global dimension, that is, 
generalized Artin--Shelter regular. However, for any such graph $\Gamma$ with at least one edge $(ij)$ 
the both these algebras are not Frobenius. 

Indeed, by definition,  an algebra $A= A_0 \oplus \dots \oplus A_n$
is Frobenius if  $\dim A_n = \dim A_0 = 1$ and the multiplication induces a non-degenerate pairing 
$A_t \otimes A_{n-t} \to A_n$ for all $t$. However, in the algebra $B^!$ we have $ B^!_{n-1} e_{ij}= 
\bigoplus k  e_{i_1}
\dots e_{i_{n-1}} e_{ij}  + \bigoplus k  e_{{i_0}{i_1}}
 e_{i_2}
\dots e_{i_{n-1}}  e_{ij}   = 0$, since in every summand at least one index appears twice, hence the pairing 
$B^!_{n-1} \otimes B^!_1 \to B^!_n$
is degenerate.  Moreover, all quadratic monomials on variables $e_i, e_{ij}$ which do not appear in the 
relations of the algebra $Q(\Gamma)$ listed in the subsection~\ref{def_QG}
are monomial relations of the algebra $Q(\Gamma)^!$ (with letters `$u$' replaced by `$e$').
In particular, we have the equalities $e_i e_{ij} = e_j e_{ij} =0$ and $e_{ij}^2=0$ in the algebra $Q(\Gamma)^!$. 
Because $B^! = \gr Q(\Gamma)^!$, every monomial basis of the vector space $B^!_{n-1}$ 
is  also a monomial basis in the vector space
$Q(\Gamma)^!_{n-1}$. Therefore,    
given two indexes $i,j$, there is a basis of the vector space $Q(\Gamma)^!_{n-1}$
consisting of the monomials 
of the form $e_{i_1}
\dots e_{i_{n-1}}$, $e_{{i_0}{i_1}}
 e_{i_2}
\dots e_{i_{n-1}}$, and $
 e_{i_1}
\dots e_{i_{n-2}} e_{ij}$, where either $i_{n-1} =i$ or $i_{n-1} =j$. 
Thus, $Q(\Gamma)^!_{n-1} e_{ij} = 0$ for every edge $(ij) \in E$, so, the algebra 
$Q(\Gamma)^!$ is not Frobenius. 
 \end{nota}




\begin{thebibliography}{Biblio}


\bibitem[BGS]{bgs}   A. Beilinson, V. Ginzburg, W. Soergel, {\it Koszul duality patterns in representation theory},  
J. Amer. Math. Soc.,  {\bf 9}  (1996),  2, p.~473--527

\bibitem[F]{fro}     R. Froberg, {\it Determination of a class of Poincare series},
        Math. Scand., {\bf 37} (1975), 1, p.~29--39

\bibitem[GRW]{grw}   I.  Gelfand, V. Retakh, R. L. Wilson,  {\it Quadratic linear algebras associated with
     factorizations of noncommutative polynomials and noncommutative differential
      polynomials},  Selecta Math. (N.S.),  {\bf  7 } (2001),  4, p.~493--523

\bibitem[GGR]{CW}     I.  Gelfand, S. Gelfand, V. Retakh,
{\it Noncommutative algebras associated to complexes and graphs},
Selecta Math (N.S.), {\bf 7}  (2001), p.~525--531

\bibitem[GGRW]{qdet}     I.  Gelfand, S. Gelfand, V. Retakh, R. Wilson,
{\it Quasideterminants}, Adv. Math., {\bf 193} (2005), p.~56--141

\bibitem[GGRSW]{5}    I. Gelfand, S. Gelfand, V. Retakh, S. Serconek, R. L. Wilson,
                 {\it Hilbert series of quadratic algebras associated with pseudo-roots of
                   noncommutative polynomials},  J. Algebra, {\bf 254}  (2002),  2, p.~279--299
		   
\bibitem[MVS]{mvs} R.  Martinez Villa, M. Saorin, {\it Koszul equivalences and dualities,}
		     Pacific J. Math., {\bf 214}  (2004),  2, p.~359--378
		   
\bibitem[N1]{nacin1} D. Nacin, {\it The Algebra $K_3$ is Koszul,}
preprint math.QA/0511589, 2005
		   
\bibitem[N2]{nacin2} D. Nacin, {\it The Algebra $P_n$ is Koszul,}
preprint math.QA/0511571, 2005

\bibitem[N3]{nacin3} D. Nacin, {\it Structural propertias of algebras arising from pseudoroots,}
Ph. D. thesis, Rutgers Univ., 2005   		   

\bibitem[Pi1]{pi1} D. Piontkovski, {\it Algebras associated to pseudo-roots of 
noncommutative polynomials are Koszul,}  
Internat. J. Algebra Comput.,  {\bf 15}  (2005),   4, p.~643--648

\bibitem[Pi2]{pi2} D. Piontkovski, 
{\it Noncommutative Gröbner bases, coherence of associative algebras, and divisibility in semigroups,} 
  Fundam. Prikl. Mat., {\bf 7}  (2001),   2, p.~495--513 [Russian]


\bibitem[PP]{pp}                   A. Polishchuk and L. Positselski,
                    {\it Quadratic algebras}, Univ. Lecture Series, {\bf 37}, AMS, 2005

\bibitem[Pr]{pri}    S. B. Priddy, {\it Koszul resolutions,}
Trans. AMS, {\bf 152} (1970), 1, p.~39--60


\bibitem[SW]{sw}   S. Serconek, R. L. Wilson, {\it The quadratic algebras associated with pseudo-roots of 
noncommutative polynomials are Koszul algebras,}   
J. Algebra,  {\bf 278}  (2004), 2, p.~473--493

\bibitem[SZ]{SZ}
         Stephenson D.R.; Zhang J.J., {\it Growth of graded
          noetherian rings,}  Proc. AMS, {\bf 125} (1997), p.~1593--1605

\bibitem[U]{ufn} V. A. Ufnarovsky, {\it Combinatorial and asymptotical
                   methods in algebra,} Sovr. probl. mat., Fund. napr.,
               {\bf 57} (1990), p.~5--177 [Russian]
               Engl. transl.: Algebra VI, Encycl. Math. Sci., Springer,
                   Berlin 1995, p.~1--196

\end{thebibliography}
\end{document}